\documentclass[12pt,twoside]{article}
\usepackage{amssymb}
\font\teneufm=eufm10 scaled \magstep1
\font\seveneufm=eufm7 scaled \magstep1
\font\fiveeufm=eufm5  scaled \magstep1
\newfam\eufmfam
\textfont\eufmfam=\teneufm
\scriptfont\eufmfam=\seveneufm
\scriptscriptfont\eufmfam=\fiveeufm
\def\frak#1{{\fam\eufmfam\relax#1}}

\newfam\msbfam
\font\tenmsb=msbm10 scaled \magstep1  \textfont\msbfam=\tenmsb
\font\sevenmsb=msbm7 scaled \magstep1 \scriptfont\msbfam=\sevenmsb
\font\fivemsb=msbm5 scaled \magstep1  \scriptscriptfont\msbfam=\fivemsb
\def\Bbb{\fam\msbfam \tenmsb}

\def\RR{{\Bbb R}}
\def\CC{{\Bbb C}}
\def\QQ{{\Bbb Q}}
\def\NN{{\Bbb N}}
\def\ZZ{{\Bbb Z}}

\def\PP{{\Bbb P}}

\def\ra{\rightarrow}

 \def\HollowBoxx #1#2#3{{\dimen0=#1 \advance\dimen0 by -#2
       \dimen1=#1 \advance\dimen1 by #3
        \vrule height 0pt depth #3 width #2
       \hskip -#3
       \vrule height #1 depth #3 width #3}}
 \def\LeftContraction{\mathord{\kern1.45pt \HollowBoxx{6pt}{3.5pt}{.4pt}}\,}

 \def\HollowBox #1#2#3{{\dimen0=#1 \advance\dimen0 by -#3
       \dimen1=#1 \advance\dimen1 by #3
        \vrule height #1 depth #3 width #3
        \vrule height 0pt depth #3 width #2
        \hskip -#3}}
 \def\RightContraction{\mathord{\, \HollowBox{6pt}{3.1pt}{.4pt}} \kern1.6pt}

\def\qed{{\hfill $\Box$}}
\newtheorem{theorem}{THEOREM}[section]

\newtheorem{lemma}[theorem]{Lemma}

\newtheorem{remark}[theorem]{Remark}
\newtheorem{proposition}[theorem]{Proposition}

\begin{document}

\begin{center}
{\Large \bf Hyperbolic Manifolds of Dimension $n$
\medskip\\
with Automorphism Group
\medskip\\ 
of Dimension $n^2-1$}
\medskip \\
\normalsize A. V. Isaev
\end{center}

\begin{quotation} \small \sl We consider complex Kobayashi-hyperbolic manifolds of dimension $n\ge 2$ for which the dimension of the group of holomorphic automorphisms is equal to $n^2-1$.  We give a complete classification of such manifolds for $n\ge 3$ and discuss several examples for $n=2$.
\end{quotation}

\thispagestyle{empty}

\pagestyle{myheadings}
\markboth{A. V. Isaev}{Automorphism Groups of Hyperbolic Manifolds}

\setcounter{section}{-1}

\section{Introduction}
\setcounter{equation}{0}

Let $M$ be a connected complex manifold and $\hbox{Aut}(M)$ the group of holomorphic automorphisms of $M$. If $M$ is Kobayashi-hyperbolic, $\hbox{Aut}(M)$ is a Lie group in the compact-open topology \cite{Ko}, \cite{Ka}. Let $d(M):=\hbox{dim}\,\hbox{Aut}(M)$. It is well-known (see \cite{Ko}, \cite{Ka}) that $d(M)\le n^2+2n$, and that $d(M)= n^2+2n$ if and only if $M$ is holomorphically equivalent to the unit ball $B^n\subset\CC^n$, where $n:=\hbox{dim}_{\CC}M$. In \cite{IKra} we studied lower automorphism group dimensions and showed that, for $n\ge 2$, there exist no hyperbolic manifolds with $n^2+3\le d(M)\le n^2+2n-1$, and that the only manifolds with $n^2<d(M)\le n^2+2$ are, up to holomorphic equivalence, $B^{n-1}\times\Delta$ (where $\Delta$ is the unit disc in $\CC$) and the 3-dimensional Siegel space (the symmetric bounded domain of type $(\hbox{III}_2)$ in $\CC^3$). Further, in \cite{I1} all manifolds with $d(M)=n^2$ were determined (for partial classifications in special cases see also \cite{GIK} and \cite{KV}). The classification in this situation is substantially richer than that for higher automorphism group dimensions.

Observe that a further decrease in $d(M)$ almost immediately leads to unclassifiable cases. For example, no good classification exists for $n=2$ and $d(M)=2$, since the automorphism group of a generic Reinhardt domain in $\CC^2$ is 2-dimensional (see also \cite{I1} for a more specific statement). While it is possible that there is some classification for $d(M)=n^2-2$, $n\ge 3$ as well as for particular pairs $d(M)$, $n$ with $d(M)< n^2-2$ (see \cite{GIK} in this regard), the case $d(M)=n^2-1$ is probably the only remaining candidate to investigate for the existence of a reasonable classification for every $n\ge 2$. It turns out that all hyperbolic manifolds with $d(M)=n^2-1$, $n\ge 2$ indeed can be explicitly described and that the case $n=2$ substantially differs from the case $n\ge 3$. In this paper we obtain a classification for $d(M)=n^2-1$, $n\ge 3$ and give examples that demonstrate some of the specifics of the case $n=2$. Our main result is the following theorem.

\begin{theorem}\label{main}\sl Let $M$ be a connected hyperbolic manifold of dimension $n\ge 3$ with $d(M)=n^2-1$. Then $M$ is holomorphically equivalent to one of the following manifolds:
\vspace{0cm}\\

\noindent (i) $B^{n-1}\times S$, where $S$ is a hyperbolic Riemann surface with $d(S)=0$;
\vspace{0cm}\\

\noindent (ii) the tube domain 
$$
\begin{array}{ll}
T:=\Bigl\{(z_1,z_2,z_3,z_4)\in\CC^4:&(\hbox{Im}\,z_1)^2+(\hbox{Im}\,z_2)^2+
\\&(\hbox{Im}\,z_3)^2-(\hbox{Im}\,z_4)^2<0,\,\hbox{Im}\,z_4>0\Bigr\}.
\end{array}
$$
(here $n=4$).
\end{theorem}

For $n=2$ in addition to the direct products specified in (i) of Theorem \ref{main} many other manifolds occur. They arise, in particular, from gluing together certain homogeneous strongly pseudoconvex real hypersurfaces in 2-dimensional complex manifolds with 3-dimensional groups of\linebreak $CR$-automorphisms. All such hypersurfaces were determined by E. Cartan \cite{C}, and our considerations for $n=2$ required an appropriate interpretation of Cartan's results (see \cite{I2}). Obtaining the classification for $n=2$ is quite lengthy, and therefore the author has decided to publish it in a separate paper. Some non-trivial examples of hyperbolic domains in $\CC^2$ and $\CC\PP^2$ with 3-dimensional automorphism groups are given in Section \ref{examples23}.

The proof of Theorem \ref{main} is organized as follows. In Section \ref{dimorbits} we determine the dimensions of the orbits of the action on $M$ of $G(M):=\hbox{Aut}(M)^c$, the connected component of the identity of $\hbox{Aut}(M)$. It turns out that, unless $M$ is homogeneous, every $G(M)$-orbit is either a real or complex hypersurface in $M$, every real hypersurface orbit is spherical and every complex hypersurface orbit is holomorphically equivalent to $B^{n-1}$ (see Proposition \ref{dim}). Note that Proposition \ref{dim} also contains some information about $G(M)$-orbits for $n=2$, in particular, it allows in this case for some real hypersurface orbits to be either Levi-flat or Levi non-degenerate non-spherical, and for some 2-dimensional orbits to be totally real rather than complex submanifolds of $M$. It turns out that such orbits indeed exist; the corresponding examples are given in Section  \ref{examples23}. 

Next, in Section \ref{sectspher} we show that real hypersurface orbits in fact cannot occur (see Proposition \ref{nospher}). First, we prove that there may be three possible kinds of such orbits and that the presence of an orbit of a particular kind determines $G(M)$ as a Lie group. Further, when we attempt to glue real hypersurface orbits together, it turns out that for any resulting hyperbolic manifold $M$, the dimension $d(M)$ is always greater than $n^2-1$. Hence all orbits are in fact complex hypersurfaces unless the manifold in question is homogeneous. Parts of the arguments in Section \ref{sectspher} apply in the case $n=2$ as well.  

In Section \ref{sectcomplex} we prove Theorem \ref{main} in the non-homogeneous case and obtain manifolds in (i) of Theorem \ref{main} (see Proposition \ref{theoremcomplex}).

In Section \ref{homog} homogeneous manifolds are considered. We show that in this case $n=4$ and obtain the tube domain in (ii) of Theorem \ref{main} (see Proposition \ref{homogeneous}). Note that Proposition \ref{homogeneous} holds for any $n\ge 2$, hence no additional homogeneous manifolds occur when $n=2$.

\section{Dimensions of Orbits}\label{dimorbits}
\setcounter{equation}{0}

The action of $G(M)=\hbox{Aut}(M)^c$ on $M$ is proper (see Satz 2.5 of \cite{Ka}), and therefore for every $p\in M$ its orbit $O(p):=\{f(p):f\in G(M)\}$ is a closed submanifold of $M$ and the isotropy subgroup $I_p:=\{f\in G(M): f(p)=p\}$ of $p$ is compact (see \cite{Ko}, \cite{Ka}). In this section we will obtain an initial classification of the $G(M)$-orbits.

Let $L_p:=\{d_pf: f\in I_p\}$ be the linear isotropy subgroup, where $d_pf$ is the differential of a map $f$ at $p$. The group $L_p$ is a compact subgroup of $GL(T_p(M),\CC)$ isomorphic to $I_p$ by means of the isotropy representation
$$
\alpha_p:\, I_p\ra L_p, \quad \alpha_p(f)=d_pf
$$
(see e.g. Satz 4.3 of \cite{Ka}). We will now prove the following proposition.

\begin{proposition}\label{dim} \sl Let $M$ be a connected hyperbolic manifold of dimension $n\ge 2$ with $d(M)=n^2-1$, and $p\in M$. Then the following holds:
\vspace{0cm}\\ 

\noindent (i) Either $M$ is homogeneous, or $O(p)$ is a real or complex closed hypersurface in $M$, or, for $n=2$, the orbit $O(p)$ is a totally real 2-dimensional closed submanifold of $M$.
\vspace{0cm}\\

\noindent (ii) If $O(p)$ is a real hypersurface, the identity component $I_p^c$ of the isotropy subgroup $I_p$ is isomorphic to $SU_{n-1}$, and $I_p$ is isomorphic to a subgroup of $\ZZ_2\times U_{n-1}$ by means of the isotropy representation $\alpha_p$. If $n\ge 3$, the orbit $O(p)$ is spherical and $I_p$ is isomorphic to a subgroup of $U_{n-1}$. If $n=2$ and $O(p)$ is strongly pseudoconvex, then it is spherical, provided $I_p$ contains more than two elements; if $n=2$ and $O(p)$ is Levi-flat, it is foliated by complex curves holomorphically equivalent to the unit disk $\Delta$.
\vspace{0cm}\\ 

\noindent (iii) If $O(p)$ is a complex hypersurface, it is holomorphically equivalent to $B^{n-1}$. If $n\ge 3$, then $I_p^c$ is isomorphic, by means of the isotropy representation $\alpha_p$, to the group $H_{k_1,k_2}^n$ of all matrices of the form 
\begin{equation}
\left(\begin{array}{cc}
a & 0\\
0 & B
\end{array}\right),\label{grouphk}
\end{equation}
where $B\in U_{n-1}$ and $a\in
(\det B)^{\frac{k_1}{k_2}}$,
for some  $k_1,k_2\in\ZZ$, $(k_1,k_2)=1$, $k_2\ne 0$. If $n=2$, then either $I_p^c$ is isomorphic, by means of the isotropy representation $\alpha_p$, to the group $H_{k_1,k_2}^2$ for some $k_1,k_2\in\ZZ$, or $L_p^c$ acts trivially on the tangent space to $O(p)$ at $p$ and $I_p^c$ is isomorphic to $U_1$ by means of the isotropy representation $\alpha_p$. If $I_p^c$ is isomorphic to $H_{k_1,k_2}^n$ for some $k_1\ne 0$, there is a real hypersurface orbit in $M$.
\vspace{0cm}\\ 

\noindent (iv) if $n=2$ and $O(p)$ is totally real, then $I_p^c$ is isomorphic to $SO_2(\RR)$ by means of the isotropy representation $\alpha_p$. 
\end{proposition}

\noindent {\bf Proof:} Let $V\subset T_p(M)$ be the tangent space to $O(p)$ at $p$. Clearly, $V$ is $L_p$-invariant. We assume now that $O(p)\ne M$ (and therefore $V\ne T_p(M)$) and consider the following three cases.
\smallskip\\

{\bf Case 1.} $d:=\hbox{dim}_{\CC}(V+iV)<n$.
\smallskip\\

Since $L_p$ is compact, one can choose coordinates in $T_p(M)$  such
 that $L_p\subset U_n$. Further, the action of $L_p$ on $T_p(M)$
is completely reducible and the subspace $V+iV$ is invariant  under this
action. Hence  $L_p$ can in fact be embedded in $U_{n-d}\times
U_d$. Since $\hbox{dim}\,O(p)\le 2d$, it follows that
$$
n^2-1\le (n-d)^2+d^2+2d,
$$
and therefore  either $d=0$ or $d=n-1$. 

If $d=0$, then $p$ is a fixed point for the action of $G(M)$ on $M$. Then $I_p=G(M)$ and $L_p$ is isomorphic to $G(M)$. Since $\hbox{dim}\,L_p=n^2-1$, we have $L_p=SU_n$. The group $SU_n$ acts transitively on directions in $T_p(M)$. Since $d(M)>0$, the manifold $M$ is non-compact. Then, by \cite{GK}, $M$ is holomorphically equivalent to $B^n$, which is clearly impossible.

Suppose that $d=n-1$. Then we have
$$
n^2-1=\hbox{dim}\,L_p+\hbox{dim}\,O(p)\le n^2-2n+2+\hbox{dim}\,O(p).
$$
Hence $\hbox{dim}\,O(p)\ge 2n-3$, that is, either $\hbox{dim}\,O(p)=2n-2$, or $\hbox{dim}\,O(p)=2n-3$. 

Suppose first that  $\hbox{dim}\,O(p)=2n-2$. In this case we have $iV=V$, hence $O(p)$ is a complex hypersurface. Then $\hbox{dim}\, L_p=(n-1)^2$. It now follows from the proof of Lemma 2.1 of \cite{IKru1} that $L_p^c$ is either $U_1\times SU_{n-1}$, or, for some $k_1$, $k_2$, the group $H_{k_1,k_2}^n$ defined in (\ref{grouphk}).
Therefore, if $n\ge 3$ or $n=2$ and $L_p^c=H_{k_1,k_2}^2$ for some $k_1$, $k_2$, then  $L_p$ acts transitively on directions in $V$, and \cite{GK} implies that $O(p)$ is holomorphically equivalent to $B^{n-1}$.

Let $n\ge 3$ and $L_p^c=U_1\times SU_{n-1}$. It then follows (see, for example, Satz 4.3 of \cite{Ka}) that $I_p':=\alpha_p^{-1}(U_1)$ is the kernel of the action of $G(M)$ on $O(p)$, in particular, $I_p'$ is normal in $G(M)$. Therefore, the factor-group $G(M)/I_p'$ acts effectively on $O(p)$. Clearly, $\hbox{dim}\,G(M)/I_p'=n^2-2$. Thus, the group $\hbox{Aut}(O(p))$ is isomorphic to $\hbox{Aut}(B^{n-1})$ (in particular, its dimension is $n^2-1$) and has a codimension 1 (possibly non-closed) subgroup. However, the Lie algebra ${\frak {su}}_{n-1,1}$ of the group $\hbox{Aut}(B^{n-1})$ does not have codimension 1 subalgebras, if $n\ge 3$ (see, e.g., \cite{EaI}). Thus, we have shown that if $n\ge 3$, then $L_p^c=H_{k_1,k_2}^n$ for some $k_1,k_2$.

Next, if $n=2$ and $L_p^c=U_1\times SU_1=U_1$, then the above argument shows that $O(p)$ is a hyperbolic 1-dimensional manifold with automorphism group of dimension at least 2. Hence $O(p)$ is holomorphically equivalent to $\Delta$ if $L_p^c=U_1$ as well.

Suppose that $I_p^c$ is isomorphic to $H_{k_1,k_2}^n$ where $k_1\ne 0$. Then $L_p^c$ acts as $U_1$ on the orthogonal complement to $V$. Therefore, in this case there are real hypersurface orbits in $M$ arbitrarily close to $O(p)$. 

Suppose now that $\hbox{dim}\,O(p)=2n-3$. In this case
$\hbox{dim}\,I_p=n^2-2n+2$. Since $L_p$ can be embedded in $U_1\times U_{n-1}$, we obtain $L_p=U_1\times U_{n-1}$. In particular, $L_p$ acts transitively on directions in $V+iV$. This is, however, impossible since $V$ is of codimension 1 in $V+iV$ and is $L_p$-invariant.   
\smallskip\\

{\bf Case 2.}  $T_p(M)=V+iV$ and $r:=\hbox{dim}_{\CC}(V\cap iV)>0$.
\smallskip\\

As above, $L_p$ can be embedded in $U_{n-r}\times U_r$ (clearly, we have
$r<n$).  Moreover,
 $V\cap iV\ne V$ and since $L_p$ preserves $V$, it follows that
$\hbox{dim}\,L_p<r^2+(n-r)^2$. We have $\hbox{dim}\,O(p)\le 2n-1$, and
therefore
$$
n^2-1<(n-r)^2+r^2+2n-1,
$$
which shows  that either $r=1$, or $r=n-1$. It then follows that 
$\hbox{dim}\,L_p<n^2-2n+2$. Therefore, we have
$$
n^2-1=\hbox{dim}\,L_p+\hbox{dim}\,O(p)<n^2-2n+2+\hbox{dim}\,O(p).
$$
Hence $\hbox{dim}\,O(p)>2n-3$. Thus, $\hbox{dim}\,O(p)=2n-1$, or
$\hbox{dim}\,O(p)=2n-2$.

Suppose that $\hbox{dim}\,O(p)=2n-1$. Let $W$ be the orthogonal complement to $V\cap iV$ in $T_p(M)$. Clearly, in this case $r=n-1$ and $\hbox{dim}_{\CC}\,W=1$. The group $L_p$ is a subgroup of $U_n$ and preserves $V$, $V\cap iV$, and $W$; hence it preserves the line $W\cap V$. Therefore, it can act only as $\pm\hbox{id}$ on $W$, that is, $L_p\subset\ZZ_2\times U_{n-1}$. Since $\hbox{dim}\,L_p=(n-1)^2-1$, we have $L_p^c=SU_{n-1}$. In particular, $L_p$ acts transitively on directions in $V\cap iV$, if $n\ge 3$. Hence, the orbit $O(p)$ is either Levi-flat or strongly pseudoconvex for all $n\ge 2$.

Suppose first that $n\ge 3$ and $O(p)$ is Levi-flat. Then $O(p)$ is foliated by connected complex manifolds. Let $M_p$ be the leaf passing through $p$. Denote by ${\frak g}$ the Lie algebra of vector fields on $O(p)$ arising from the action of $G(M)$, and let ${\frak l}_p\subset{\frak g}$ be the subspace consisting of all vector fields tangent to $M_p$ at $p$. Since vector fields in ${\frak l}_p$ remain tangent to $M_p$ at each point in $M_p$, the subspace ${\frak l}_p$ is in fact a Lie subalgebra of ${\frak g}$. It follows from the definition of ${\frak l}_p$ that $\hbox{dim}\,{\frak l}_p=n^2-2$. Denote by $H_p$ the (possibly non-closed) connected subgroup of $G(M)$ with Lie algebra ${\frak l}_p$. It is straightforward to verify that the group $H_p$ acts on $M_p$ by holomorphic transformations and that $I_p^c\subset H_p$. If some non-trivial element $g\in H_p$ acts trivially on $M_p$, then $g\in I_p$, and corresponds to the non-trivial element in $\ZZ_2$ (recall that $L_p\subset \ZZ_2\times U_{n-1}$). Thus, either $H_p$ or $H_p/\ZZ_2$ acts effectively on $M_p$ (the former case occurs if $g_p\not\in H_p$, the latter if $g_p\in H_p$). The group $L_p$ acts transitively on directions in the tangent space $V\cap iV$ to $M_p$, and it follows from \cite{GK} that $M_p$ is holomorphically equivalent to $B^{n-1}$. Therefore, the group $\hbox{Aut}(M_p)$ is isomorphic to $\hbox{Aut}(B^{n-1})$ (in particular, its dimension is $n^2-1$) and has a codimension 1 (possibly non-closed) subgroup. 
However, as we noted above, the Lie algebra of $\hbox{Aut}(B^{n-1})$ does not have codimension 1 subalgebras, if $n\ge 3$. Thus, $O(p)$ is strongly pseudoconvex. Hence, $L_p$ acts trivially on $W$ and therefore $L_p\subset U_{n-1}$. Since $L_p^c=SU_{n-1}$, the dimension of the stability group of $O(p)$ at $p$ is greater than or equal to $(n-1)^2-1$, which for $n\ge 3$ implies that $p$ is an umbilic point of $O(p)$ (see e.g. \cite{EzhI}). The homogeneity of $O(p)$ now yields that $O(p)$ is spherical, if $n\ge 3$. For $n=2$ the above argument shows that $O(p)$ is foliated by connected hyperbolic complex curves  with automorphism group of dimension at least 2, that is, by complex curves holomorphically equivalent to $\Delta$.

If $n=2$, the orbit $O(p)$ is Levi non-degenerate and $I_p$ contains more than two elements, then arguing as in the proof of Lemma 3.3 of \cite{IKru2}, we obtain that $O(p)$ is spherical. Alternatively, this fact can be derived from the classification in \cite{C}.   

Suppose now that $\hbox{dim}\,O(p)=2n-2$. Since $T_p(M)=V+iV$, the orbit $O(p)$ is not a complex hypersurface. Therefore, $r=n-2$, which is only possible for $n=3$ (recall that we have either $r=1$, or $r=n-1$). In this case $\hbox{dim}\,L_p=4$ and, arguing as in the proof of Lemma 2.1 of \cite{IKru1}, we see that $L_p$ acts transitively on directions in the orthogonal complement $W$ to $V\cap iV$ in $T_p(M)$. This is, however, impossible since $L_p$ must preserve $W\cap V$.
\smallskip\\

{\bf Case 3.}  $T_p(M)=V\oplus iV$.
\smallskip\\

In this case  $\hbox{dim}\, V=n$ and $L_p$ can be embedded in the real orthogonal group $O_n(\RR)$,
and therefore
$$
\hbox{dim}\,L_p+\hbox{dim}\, O(p)\le \frac{n(n-1)}{2}+n.
$$
Hence, for $n\ge 3$, we have $\hbox{dim}\,L_p+\hbox{dim}\, O(p)<n^2-1$
which is impossible.

Assume now that $n=2$. If $\hbox{dim}\,L_p=0$, we get a contradiction as above. Hence $\hbox{dim}\,L_p=1$ and $L_p^c=SO_2(\RR)$.

The proof of the proposition  is complete.\qed
\smallskip\\

\section{Real Hypersurface Orbits}\label{sectspher}
\setcounter{equation}{0}

In this section we will deal with real hypersurface orbits and eventually show that they do not occur. Our goal is to prove the following proposition.

\begin{proposition}\label{nospher}\sl Let $M$ be a connected hyperbolic manifold of dimension $n\ge 3$ with $d(M)=n^2-1$. Then no orbit in $M$ is a real hypersurface.
\end{proposition}

\noindent {\bf Proof:} Recall that every real hypersurface orbit is spherical. First, we narrow down the class of all possible spherical orbits.

\begin{lemma}\label{spherorbitsprop}\sl Let $M$ be a connected hyperbolic manifold of dimension $n\ge 3$ with $d(M)=n^2-1$. Assume that for a point $p\in M$ its orbit $O(p)$ is spherical. Then $O(p)$ is $CR$-equivalent to one of the following hypersurfaces:
\begin{equation}
\begin{array}{ll}
\hbox{(i)}& \hbox{a lens manifold ${\cal L}_m:=S^{2n-1}/\ZZ_m$ for some $m\in\NN$},\\
\hbox{(ii)} & \sigma:=\left\{(z',z_n)\in\CC^{n-1}\times\CC:\hbox{Re}\,z_n=|z'|^2\right\},\\
\hbox{(iii)} & \delta:=\hbox{$\left\{(z',z_n)\in\CC^{n-1}\times\CC: |z_n|=\exp\left(|z'|^2\right)\right\}$},\\
\hbox{(iv)} & \omega:=\left\{(z',z_n)\in\CC^{n-1}\times\CC:|z'|^2+\exp\left(\hbox{Re}\,z_n\right)=1\right\},\\
\hbox{(v)} & \varepsilon_{\alpha}:=\hbox{$\left\{(z',z_n)\in\CC^{n-1}\times\CC: |z'|^2+|z_n|^{\alpha}=1,\, z_n\ne 0\right\}$,}\\
&\hbox{for some $\alpha>0$.}
\end{array}\label{classificationspherorb}
\end{equation}
\end{lemma}

\noindent {\bf Proof of Lemma \ref{spherorbitsprop}:} The proof is similar to that of Proposition 3.1 of \cite{I1}. For a connected Levi non-degenerate $CR$-manifold $Q$ denote by $\hbox{Aut}_{CR}(Q)$ the Lie group of its $CR$-automorphisms. Let $\tilde O(p)$ be the universal cover of $O(p)$. The connected component of the identity $\hbox{Aut}_{CR}(O(p))^c$ of $\hbox{Aut}_{CR}(O(p))$ acts transitively on $O(p)$ and therefore its universal cover $\widetilde{\hbox{Aut}}_{CR}(O(p))^c$ acts transitively on $\tilde O(p)$. Let $G$ be the (possibly non-closed) subgroup of $\hbox{Aut}_{CR}(\tilde O(p))$ that consists of all $CR$-automorphisms of $\tilde O(p)$ generated by this action. Observe that $G$ is a Lie group isomorphic to the factor-group of $\widetilde{\hbox{Aut}}_{CR}(O(p))^c$ by a discrete central subgroup. Let $\Gamma\subset \hbox{Aut}_{CR}(\tilde O(p))$ be the discrete subgroup whose orbits are the fibers of the covering $\tilde O(p)\ra O(p)$. The group $\Gamma$ acts freely properly discontinuously on $\tilde O(p)$, lies in the centralizer of $G$ in $\hbox{Aut}_{CR}(\tilde O(p))$ and is isomorphic to $H/H^c$, with $H=\pi^{-1}(I_p)$, where $\pi:\widetilde{\hbox{Aut}}_{CR}(O(p))^c\ra\hbox{Aut}_{CR}(O(p))^c$ is the covering map.   

The manifold $\tilde O(p)$ is spherical, and there is a local $CR$-isomorphism $\Pi$ from $\tilde O(p)$ onto a domain $D\subset S^{2n-1}$. By Proposition 1.4 of \cite{BS}, $\Pi$ is a covering map. Further, for every $f\in\hbox{Aut}_{CR}(\tilde O(p))$ there is $g\in\hbox{Aut}(D)$ such that 
\begin{equation}
g\circ \Pi=\Pi\circ f.\label{liftspher}
\end{equation}
Since $\tilde O(p)$ is homogeneous, (\ref{liftspher}) implies that $D$ is homogeneous as well, and $\hbox{dim}\,\hbox{Aut}_{CR}(\tilde O(p))=\hbox{dim}\,\hbox{Aut}_{CR}(D)$.

Clearly, $\dim\hbox{Aut}_{CR}(O(p))\ge n^2-1$ and therefore we have $\dim\hbox{Aut}_{CR}(D)\ge n^2-1$. All homogeneous domains in $S^{2n-1}$ are listed in Theorem 3.1 in \cite{BS}. It is not difficult to exclude from this list all the domains with automorphism group of dimension less than $n^2-1$. This gives that $D$ is $CR$-equivalent to one of the following domains:
$$
\begin{array}{ll}
\hbox{(a)}& S^{2n-1},\\
\hbox{(b)}& S^{2n-1}\setminus\{\hbox{point}\},\\
\hbox{(c)}& S^{2n-1}\setminus\{z_n=0\}.
\end{array}
$$
Thus, $\tilde O(p)$ is respectively one of the following manifolds:
$$
\begin{array}{ll}
\hbox{(a)}& S^{2n-1},\\
\hbox{(b)}& \sigma,\\
\hbox{(c)}& \omega.
\end{array}
$$

If $\tilde O(p)=S^{2n-1}$, then by Proposition 5.1 of \cite{BS} the orbit $O(p)$ is $CR$-equivalent to a lens manifold as in (i) of (\ref{classificationspherorb}). 

Suppose next that $\tilde O(p)=\sigma$. The group $\hbox{Aut}_{CR}(\sigma)$ consists of all maps of the form
\begin{equation}
\begin{array}{lll}
z' & \mapsto & \lambda Uz'+a,\\
z_n & \mapsto & \lambda^2z_n+2\lambda\langle Uz',a\rangle+|a|^2+i\alpha,
\end{array}\label{thegroupsphpt}
\end{equation}
where $U\in U_{n-1}$, $a\in\CC^{n-1}$, $\lambda\in\RR^*$, $\alpha\in\RR$, and $\langle\cdot\,,\cdot\rangle$ is the inner product in $\CC^{n-1}$. It then follows that $\hbox{Aut}_{CR}(\sigma)=CU_{n-1}\ltimes N$, where $CU_{n-1}$ consists of all maps of the form (\ref{thegroupsphpt}) with $a=0$, $\alpha=0$, and $N$ is the Heisenberg group consisting of the maps of the form (\ref{thegroupsphpt}) with $U=\hbox{id}$ and $\lambda=1$.

Further, description (\ref{thegroupsphpt}) implies that $\hbox{dim}\,\hbox{Aut}_{CR}(\sigma)=n^2+1$, and therefore $n^2-1\le \hbox{dim}\,G\le n^2+1$. If $\hbox{dim}\,G=n^2+1$, then we have $G=\hbox{Aut}_{CR}(\sigma)^c$, and hence $\Gamma$ is a central subgroup of $\hbox{Aut}_{CR}(\sigma)^c$. Since the center of $\hbox{Aut}_{CR}(\sigma)^c$ is trivial, so is $\Gamma$. Thus, in this case $O(p)$ is $CR$-equivalent to the hypersurface $\sigma$.

Assume now that $n^2-1\le\hbox{dim}\,G\le n^2$. Since $G$ acts transitively on $\sigma$, we have $N\subset G$. Furthermore, since $G$ is of codimension 1 or 2 in $\hbox{Aut}_{CR}(\sigma)$, it either contains the subgroup $SU_{n-1}\ltimes N$, or $n=3$ and $G$ contains a subgroup of the form $L\ltimes N$, where $L$ is conjugate to $U_1\times U_1$ in $U_2$. By Proposition 5.6 of \cite{BS}, we have $\Gamma\subset U_{n-1}\ltimes N$. The centralizer of $SU_{n-1}\ltimes N$ in $U_{n-1}\ltimes N$ and that of $L\ltimes N$ in $U_2\ltimes N$ consist of all maps of the form
\begin{equation}
\begin{array}{lll}
z' &\mapsto & z',\\
z_n &\mapsto & z_n+i\alpha,
\end{array}\label{center}
\end{equation}
where $\alpha\in\RR$. Since $\Gamma$ acts freely properly discontinuously on $\sigma$, it is generated by a single map of the form (\ref{center}) with $\alpha=\alpha_0\in\RR^*$. The hypersurface $\sigma$ covers the hypersurface
\begin{equation}
\left\{(z',z_n)\in\CC^{n-1}\times\CC: |z_n|=\exp\left(\frac{2\pi}{\alpha_0}|z'|^2\right)\right\}\label{intermediate1}
\end{equation}
by means of the map
\begin{equation}
\begin{array}{lll}
z' & \mapsto & z',\\
z_n & \mapsto & \exp\left(\displaystyle\frac{2\pi}{\alpha_0} z_n\right),
\end{array}\label{coverrr}
\end{equation}
and the fibers of this map are the orbits of $\Gamma$. Hence $O(p)$ is $CR$-equivalent to hypersurface (\ref{intermediate1}). Replacing if necessary $z_n$ by $1/z_n$ we obtain that $O(p)$ is $CR$-equivalent to  the hypersurface $\delta$.

Suppose finally that $\tilde O(p)=\omega$. First, we will determine the group $\hbox{Aut}_{CR}(\omega)$. The general form of a $CR$-automorphism of $S^{2n-1}\setminus\{z_n=0\}$ is given by the formula 
$$
\begin{array}{lll}
z'&\mapsto&\displaystyle\frac{Az'+b}{cz'+d},\\
\vspace{0mm}&&\\
z_n&\mapsto&\displaystyle
\frac{e^{i\beta}z_n}{cz'+d},
\end{array}
$$
where
$$
\left(\begin{array}{cc}
A& b\\
c& d
\end{array}
\right)
\in SU_{n-1,1},\quad\beta\in\RR,
$$
and the covering map $\Pi$ by the formula      
$$
\begin{array}{l}
z'\mapsto z',\\
z_n\mapsto \exp\left(\displaystyle\frac{z_n}{2}\right).
\end{array}
$$
Using (\ref{liftspher}) we then obtain the general form of a $CR$-automorphism of $\omega$ as follows
\begin{equation}
\begin{array}{lll}
z'&\mapsto&\displaystyle\frac{Az'+b}{cz'+d},\\
\vspace{0mm}&&\\
z_n&\mapsto&\displaystyle z_n-2\ln(cz'+d)+i\beta,
\end{array}\label{autgrpcov}
\end{equation}
where
$$
\left(\begin{array}{cc}
A& b\\
c& d
\end{array}
\right)
\in SU_{n-1,1},\quad\beta\in\RR.
$$
In particular, $\hbox{Aut}_{CR}(\omega)$ is a connected group of dimension $n^2$, and therefore $n^2-1\le\hbox{dim}\,G\le n^2$.

Thus, either $G=\hbox{Aut}_{CR}(\omega)$, or $G$ coincides with the subgroup of $\hbox{Aut}_{CR}(\omega)$ given by the condition $\beta=0$ in formula (\ref{autgrpcov}). In either case, the centralizer of $G$ in $\hbox{Aut}_{CR}(\omega)$ consists of all maps of the form (\ref{center}). Hence $\Gamma$ is generated by a single such map with $\alpha=\alpha_0\in\RR$. If $\alpha_0=0$, the orbit $O(p)$ is $CR$-equivalent to $\omega$. 
Let $\alpha_0\ne 0$. The hypersurface $\omega$ covers the hypersurface
\begin{equation}
\left\{(z',z_n)\in\CC^{n-1}\times\CC: |z'|^2+|z_n|^{\frac{\alpha_0}{2\pi}}=1,\, z_n\ne 0\right\}\label{intermediate2}
\end{equation}
by means of map (\ref{coverrr}). Since the fibers of this map are the orbits of $\Gamma$, it follows that $O(p)$ is $CR$-equivalent to hypersurface (\ref{intermediate2}). Replacing if necessary $z_n$ by $1/z_n$, we obtain that $O(p)$ is $CR$-equivalent to the hypersurface 
$\varepsilon_{\alpha}$ for some $\alpha>0$. 

The proof of Lemma \ref{spherorbitsprop} is complete.\qed 
\smallskip\\

\begin{remark}\label{spher2}\rm For $n=2$ there is an additional possibility for $D$ that has to be taken into the account. Namely, $S^3\setminus \RR^2$ has a 3-dimensional automorphism group arising from the natural transitive action of $O^c_{2,1}(\RR)$ by fractional-linear transformations (see Section \ref{examples23}).  
\end{remark}

We will now show that in most cases the presence of a spherical orbit of a particular kind in $M$ determines the group $G(M)$ as a Lie group. Suppose that for some $p\in M$ the orbit $O(p)$ is spherical, and let ${\frak m}$ be the manifold from list (\ref{classificationspherorb}) to which $O(p)$ is $CR$-equivalent (we say that ${\frak m}$ is the {\it model}\, for $O(p)$). Since $G(M)$ acts effectively on $O(p)$, the $CR$-equivalence induces an isomorphism between $G(M)$ and a (possibly non-closed) connected $(n^2-1)$-dimensional subgroup $R_{\frak m}$ of $\hbox{Aut}_{CR}({\frak m})$ (this subgroup a priori depends on the choice of the $CR$-equivalence). 

We need the following lemma.

\begin{lemma}\label{groupsdeterm} \sl${}$\linebreak

\noindent(i) $R_{S^{2n-1}}$ is conjugate to $SU_n$ in $\hbox{Aut}(B^n)$,  and $R_{{\cal L}_m}=SU_n/(SU_n\cap\ZZ_m)$ for $m>1$; 

\noindent (ii) $R_{\sigma}=SU_{n-1}\ltimes N$;

\noindent (iii) $R_{\delta}$ consists of all maps of the form   
$$
\begin{array}{lll}
z' & \mapsto & Uz'+a,\\
z_n & \mapsto &e^{i\beta}\exp\Bigl(2\langle Uz',a\rangle+|a|^2\Bigr)z_n,
\end{array}
$$
where $U\in SU_{n-1}$, $a\in\CC^{n-1}$, $\beta\in\RR$;

\noindent (iv) $R_{\omega}$ consists of all maps of the form (\ref{autgrpcov}) with $\beta=0$;

\noindent (v) $R_{\varepsilon_{\alpha}}$ consists of all maps of the form
\begin{equation}
\begin{array}{lll}
z'&\mapsto&\displaystyle\frac{Az'+b}{cz'+d},\\
\vspace{0mm}&&\\
z_n&\mapsto&\displaystyle
\frac{z_n}{(cz'+d)^{2/\alpha}},
\end{array}\label{repsilon}
\end{equation}
where
$$
\left(\begin{array}{cc}
A& b\\
c& d
\end{array}
\right)
\in SU_{n-1,1}.
$$ 
\end{lemma}

\noindent{\bf Proof of Lemma \ref{groupsdeterm}:} Suppose first that ${\frak m}={\cal L}_m$, for some $m\in\NN$. Then $O(p)$ is compact and, since $I_p$ is compact as well, it follows that $G(M)$ is compact. Assume first that $m=1$. In this case $R_{S^{2n-1}}$ is a subgroup of $\hbox{Aut}_{CR}(S^{2n-1})=\hbox{Aut}(B^n)$.
Since $R_{S^{2n-1}}$ is compact, it is conjugate to a subgroup of $U_n$, which is a maximal compact subgroup in $\hbox{Aut}(B^n)$. Since both $R_{S^{2n-1}}$ is $(n^2-1)$-dimensional, it is conjugate to $SU_n$. Suppose now that $m>1$. It is straightforward to determine the group $\hbox{Aut}_{CR}\left({\cal L}_m\right)$ by lifting $CR$-automorphisms of ${\cal L}_m$ to its universal cover $S^{2n-1}$. This group is $U_n/\ZZ_m$ acting on $\CC^n\setminus\{0\}/\ZZ_m$ in the standard way. Since $R_{{\cal L}_m}$ is of codimension 1 in $\hbox{Aut}_{CR}\left({\cal L}_m\right)$, we obtain $R_{{\cal L}_m}=SU_n/(SU_n\cap\ZZ_m)$.     

Assume now that ${\frak m}=\sigma$. The group $\hbox{Aut}_{CR}(\sigma)$ consists of all maps of the form (\ref{thegroupsphpt}) and has dimension $n^2+1$. 
Since $R_{\sigma}$ acts transitively on $\sigma$, it contains the subgroup $N$ (see the proof of Proposition \ref{spherorbitsprop}). 
Furthermore, $R_{\sigma}$ is a codimension 2 subgroup of $\hbox{Aut}_{CR}(\sigma)$, and thus either is the group $SU_{n-1}\ltimes N$, or, for $n=3$, $R_{\sigma}\cap (U_2\ltimes N)=L\ltimes N$, where $L$ is conjugate to $U_1\times U_1$ in $U_2$. By (ii) of Proposition \ref{dim}, $I_p^c$ is isomorphic to $SU_{n-1}$, hence the latter case in fact does not occur.

Next, the group $\hbox{Aut}_{CR}(\delta)$ can be determined by considering the universal cover of $\delta$ (see the proof of Proposition \ref{spherorbitsprop}) and consists of all maps of the form 
\begin{equation}
\begin{array}{lll}
z' & \mapsto & Uz'+a,\\
z_n & \mapsto &e^{i\beta}\exp\Bigl(2\langle Uz',a\rangle+|a|^2\Bigr)z_n,
\end{array}\label{gdelta}
\end{equation}
where $U\in U_{n-1}$, $a\in\CC^{n-1}$, $\beta\in\RR$. This group has dimension $n^2$, and hence $R_{\delta}$ is of codimension 1 in $\hbox{Aut}_{CR}(\delta)$. Since $R_{\delta}$ acts transitively on $\delta$, it consists of all maps of the form (\ref{gdelta}) with $U\in SU_{n-1}$.

Assume now that ${\frak m}=\omega$. The only codimension 1 subgroup of $\hbox{Aut}_{CR}(\delta)$ is given by maps with $\beta=0$ in formula (\ref{autgrpcov}).

Let finally ${\frak m}=\varepsilon_{\alpha}$. The group $\hbox{Aut}_{CR}(\varepsilon_{\alpha})$ consists of all maps of the form 
\begin{equation}
\begin{array}{lll}
z'&\mapsto&\displaystyle\frac{Az'+b}{cz'+d},\\
\vspace{0mm}&&\\
z_n&\mapsto&\displaystyle
\frac{e^{i\beta}z_n}{(cz'+d)^{2/\alpha}},
\end{array}\label{autvarepsilon}
\end{equation}
where
$$
\left(\begin{array}{cc}
A& b\\
c& d
\end{array}
\right)
\in SU_{n-1,1},\quad\beta\in\RR,
$$
and its only codimension 1 subgroup is given by $\beta=0$.
    
The proof of Lemma \ref{groupsdeterm} is complete.\qed

We will now finish the proof of Proposition \ref{nospher}. Our argument is similar to that in Section 4 of \cite{I1}. For completeness of our exposition, we will repeat it here in detail.

Suppose that for some $p\in M$ the orbit $O(p)$ is $CR$-equivalent to a lens manifold ${\cal L}_m$. In this case $G(M)$ is compact, hence there are no complex hypersurface orbits and the model for every orbit is a lens manifold. Assume first that $m=1$. Then $M$ admits an effective action of $SU_n$ by holomorphic transformations and therefore is holomorphically equivalent to one of the manifolds listed in \cite{IKru2}. However, none of the manifolds on the list in \cite{IKru2} with $n\ge 3$ is hyperbolic and has $(n^2-1)$-dimensional automorphism group.

Assume now that $m>1$. Let $f: O(p)\ra {\cal L}_m$ be a $CR$-isomorphism. Then we have 
\begin{equation}
f(gq)=\varphi(g)f(q),\label{equivar}
\end{equation}
where $q\in O(p)$, for some Lie group isomorphism $\varphi: G(M)\ra SU_n/(SU_n\cap\ZZ_m)$. The $CR$-isomorphism $f$ extends to a biholomorphic map from a neighborhood $U$ of $O(p)$ in $M$ onto a neighborhood $W$ of ${\cal L}_m$ in $\CC^n\setminus\{0\}/\ZZ_m$. Since $G(M)$ is compact, one can choose $U$ to be a connected union of $G(M)$-orbits. Then property (\ref{equivar}) holds for the extended map, and therefore every $G(M)$-orbit in $U$ is taken onto an $SU_n/(SU_n\cap\ZZ_m)$-orbit in $\CC^n\setminus\{0\}/\ZZ_m$ by this map. Thus, $W=S_r^R/\ZZ_m$ for some $0\le r<R<\infty$, where $S_r^R:=\left\{z\in\CC^n: r<|z|<R\right\}$ is a spherical shell.

Let $D$ be a maximal domain in $M$ such that there exists a biholomorphic map $f$ from $D$ onto $S_r^R/\ZZ_m$ for some $r,R$, satisfying (\ref{equivar}) for all $g\in G(M)$ and $q\in D$. As was shown above, such a domain $D$ exists. Assume that $D\ne M$ and let $x$ be a boundary point of $D$. Consider the orbit $O(x)$. Let ${\cal L}_k$ for some $k>1$ be the model for $O(x)$ and $f_1:O(x)\ra{\cal L}_k$ a $CR$-isomorphism satisfying (\ref{equivar}) for $g\in G(M)$, $q\in O(x)$ and an isomorphism $\varphi_1:G(M)\ra SU_n/(SU_n\cap\ZZ_k)$ in place of $\varphi$. The map $f_1$ can be holomorphically extended to a neighborhood $V$ of $O(x)$ that one can choose to be a connected union of $G(M)$-orbits. The extended map satisfies (\ref{equivar}) for $g\in G(M)$, $q\in V$ and $\varphi_1$ in place of $\varphi$. For $s\in V\cap D$ we consider the orbit $O(s)$. The maps $f$ and $f_1$ take $O(s)$ into some surfaces $r_1S^{2n-1}/\ZZ_m$ and $r_2S^{2n-1}/\ZZ_k$, respectively, with $r_1,r_2>0$.
Hence $F:=f_1\circ f^{-1}$ maps $r_1S^{2n-1}/\ZZ_m$ onto $r_2S^{2n-1}/\ZZ_k$. Since ${\cal L}_m$ and ${\cal L}_k$ are not $CR$-equivalent for distinct $m$, $k$, we obtain $k=m$. Furthermore, every $CR$-isomorphism between $r_1S^{2n-1}/\ZZ_m$ and $r_2S^{2n-1}/\ZZ_m$ has the form $[z]\mapsto [r_2/r_1Uz]$, where $U\in U_n$, and $[z]\in\CC^n\setminus\{0\}/\ZZ_m$ denotes the equivalence class of a point $z\in\CC^n\setminus\{0\}$. Therefore, $F$ extends to a holomorphic automorphism of $\CC^n\setminus\{0\}/\ZZ_m$.

We claim that $V$ can be chosen so that $D\cap V$ is connected and\linebreak $V\setminus(D\cup O(x))\ne\emptyset$. Indeed, since $O(x)$ is strongly pseudoconvex and closed in $M$, for $V$ small enough we have $V=V_1\cup V_2\cup O(x)$, where $V_j$ are open connected non-intersecting sets. For each $j$, $D\cap V_j$ is a union of $G(M)$-orbits and therefore is mapped by $f$ onto a union of the quotients of some spherical shells. If there are more than one such factored shells, then there is a factored shell such that the closure of its inverse image under $f$ is disjoint from $O(x)$, and hence $D$ is disconnected which contradicts the definition of $D$. Thus, $D\cap V_j$ is connected for $j=1,2$, and, if $V$ is sufficiently small, then each $V_j$ is either a subset of $D$ or is disjoint from it. If $V_j\subset D$ for $j=1,2$, then $M=D\cup V$ is compact, which is impossible since $M$ is hyperbolic and $d(M)>0$. Therefore, for some $V$ there is only one $j$ for which $D\cap V_j\ne\emptyset$. Thus, $D\cap V$ is connected and $V\setminus(D\cup O(x))\ne\emptyset$, as required. 

Setting now
\begin{equation}
\tilde f:=\Biggl\{\begin{array}{l}
f\hspace{1.7cm}\hbox{on $D$}\\
F^{-1} \circ f_1\hspace{0.45cm}\hbox{on $V$},
\end{array}\label{extens}
\end{equation}
we obtain a biholomorphic extension of $f$ to $D\cup V$. By construction, $\tilde f$ satisfies (\ref{equivar}) for $g\in G(M)$ and $q\in D\cup V$. Since $D\cup V$ is strictly larger than $D$, we obtain a contradiction with the maximality of $D$. Thus, we have shown that in fact $D=M$, and hence $M$ is holomorphically equivalent to $S_r^R/\ZZ_m$. However, in this case $d(M)=n^2$, which is impossible.

The orbit gluing procedure utilized above can in fact be applied in a very general setting. We will now describe it in full generality (see also \cite{I1}), assuming that every orbit in $M$ is a real hypersurface. The procedure comprises the following steps:
\vspace{0cm}\\

\noindent (1). Start with a real hypersurface orbit $O(p)$ with model ${\frak m}$ and consider a real-analytic $CR$-isomorphism $f:O(p)\ra{\frak m}$ that satisfies (\ref{equivar}) for all $g\in G(M)$ and $q\in O(p)$, where $\varphi:G(M)\ra R_{\frak m}$ is a Lie group isomorphism. 
\vspace{0cm}\\

\noindent (2). Verify that for every model ${\frak m}'$ the group $R_{{\frak m}'}$ acts by holomorphic transformations with real hypersurface orbits on a domain ${\cal D}\subset\CC^n$ that contains ${\frak m}'$ and that every orbit of the action is $CR$-equivalent to ${\frak m}'$.
\vspace{0cm}\\

\noindent (3). Observe that $f$ can be extended to a biholomorphic map from a $G(M)$-invariant connected neighborhood of $O(p)$ in $M$ onto an $R_{\frak m}$-invariant neighborhood of ${\frak m}$ in ${\cal D}$. First of all, extend $f$ to some neighborhood $U$ of $O(p)$ to a biholomorphic map onto a neighborhood $W$ of ${\frak m}$ in $\CC^n$. Let $W'=W\cap {\cal D}$ and $U'=f^{-1}(W')$. Fix $s\in U'$ and $s_0\in O(s)$. Choose $h_0\in G(M)$ such that $s_0=h_0s$ and define $f(s_0):=\varphi(h_0)f(s)$. To see that $f$ is well-defined at $s_0$, suppose that for some $h_1\in G(M)$, $h_1\ne h_0$, we have $s_0=h_1s$, and show that $\varphi(h)$ fixes $f(s)$, where $h:=h_1^{-1}h_0$. Indeed, for every $g\in G(M)$ identity (\ref{equivar}) holds for $q\in U_g$, where $U_g$ is the connected component of $g^{-1}(U')\cap U'$ containing $O(p)$. Since $h\in I_s$, we have $s\in U_h$ and the application of (\ref{equivar}) to $h$ and $s$ yields that $\varphi(h)$ fixes $f(s)$, as required. Thus, $f$ extends to $U'':=\cup_{q\in U'}O(q)$. The extended map satisfies (\ref{equivar}) for all $g\in G(M)$ and $q\in U''$.
\vspace{0cm}\\  

\noindent (4). Consider a maximal $G(M)$-invariant domain $D\subset M$ from which there exists a biholomorphic map $f$ onto an $R_{\frak m}$-invariant domain in ${\cal D}$ satisfying (\ref{equivar}) for all $g\in G(M)$ and $q\in D$. The existence of such a domain is guaranteed by the previous step. Assume that $D\ne M$ and consider $x\in\partial D$. Let ${\frak m}_1$ be the model for $O(x)$ and let $f_1:O(x)\ra{\frak m}_1$ be a real-analytic $CR$-isomorphism satisfying (\ref{equivar}) for all $g\in G(M)$, $q\in O(x)$ and some Lie group isomorphism $\varphi_1: G(M)\ra R_{{\frak m}_1}$ in place of $\varphi$. 
Let ${\cal D}_1$ be the domain in $\CC^n$ containing ${\frak m}_1$ on which $R_{{\frak m}_1}$ acts by holomorphic transformations with real hypersurface orbits $CR$-equivalent to ${\frak m}_1$. As in (3), extend $f_1$ to a biholomorphic map from a connected $G(M)$-invariant neighborhood $V$ of $O(x)$ onto an $R_{{\frak m}_1}$-invariant neighborhood of ${\frak m}_1$ in ${\cal D}_1$. The extended map satisfies (\ref{equivar}) for all $g\in G(M)$, $q\in V$ and $\varphi_1$ in place of $\varphi$. Consider $s\in V\cap D$. The maps $f$ and $f_1$ take $O(s)$ onto an $R_{\frak m}$-orbit in ${\cal D}$ and an $R_{{\frak m}_1}$-orbit in ${\cal D}_1$, respectively. Then $F:=f_1\circ f^{-1}$ maps the $R_{\frak m}$-orbit onto the $R_{{\frak m}_1}$-orbit. Since all models are pairwise $CR$ non-equivalent, we obtain ${\frak m}_1={\frak m}$. 
\vspace{0cm}\\   

\noindent (5). Show that $F$ extends to a holomorphic automorphism of ${\cal D}$.
\vspace{0cm}\\ 

\noindent (6). Show that $V$ can be chosen so that $D\cap V$ is connected and\linebreak $V\setminus(D\cup O(x))\ne\emptyset$. This follows from the hyperbolicity of $M$ and the existence of a neighborhood $V'$ of $O(x)$ such that $V'=V_1\cup V_2\cup O(x)$, where $V_j$ are open connected non-intersecting sets. The existence of such $V'$ follows from the strong pseudoconvexity of ${\frak m}$.
\vspace{0cm}\\ 

\noindent (7). Use formula (\ref{extens}) to extend $f$ to $D\cup V$ thus obtaining a contradiction with the maximality of $D$. This shows that in fact $D=M$ and hence $M$ is biholomorphically equivalent to an $R_{\frak m}$-invariant domain in ${\cal D}$. In all the cases below the determination of $R_{\frak m}$-invariant domains will be straightforward.
\vspace{0cm}\\ 

Assume first that every orbit in $M$ is a real hypersurface. Let first ${\frak m}=\sigma$. Clearly, the group $R_{\sigma}$ acts with real hypersurface orbits on all of $\CC^n$, so in this case ${\cal D}=\CC^n$. The $R_{\sigma}$-orbit of every point in $\CC^n$ is of the form
$$
\left\{(z',z_n)\in\CC^{n-1}\times\CC:\hbox{Re}\,z_n=|z'|^2+r\right\},
$$
where $r\in\RR$, and every $R_{\sigma}$-invariant domain in $\CC^n$ is given by
$$
{\frak S}_r^R:=\left\{(z',z_n)\in\CC^{n-1}\times\CC: r+|z'|^2<\hbox{Re}\,z_n<R+|z'|^2\right\},
$$
where $-\infty\le r<R\le\infty$. Every $CR$-isomorphism between two $R_{\sigma}$-orbits is a composition of a map of the form (\ref{thegroupsphpt}) and a translation in the $z_n$-variable. Therefore, $F$ in this case extends to a holomorphic automorphism of $\CC^n$. Now our gluing procedure implies that $M$ is holomorphically equivalent to ${\frak S}_r^R$ for some $-\infty\le r<R\le\infty$. Therefore, $M$ is holomorphically equivalent either to the domain 
$$
{\frak S}:=\Bigl\{(z',z_n)\in\CC^{n-1}\times\CC: -1+|z'|^2<\hbox{Re}\,z_n<|z'|^2\Bigr\},
$$
or (for $R=\infty$) to $B^n$. The latter is clearly impossible; the former is impossible either since $d({\frak S})=n^2$ (see e.g. \cite{I1}).

Assume next that ${\frak m}=\delta$. Again, we have ${\cal D}=\CC^n$. The $R_{\delta}$-orbit of every point in $\CC^n$ has the form 
$$
\left\{(z',z_n)\in\CC^{n-1}\times\CC: |z_n|=r\exp\left(|z'|^2\right)\right\},
$$
where $r>0$, and hence every $R_{\delta}$-invariant domain in $\CC^n$ is given by
$$
D_r^R:=\Bigl\{(z',z_n)\in\CC^{n-1}\times\CC: r\exp\left({|z'|^2}\right)<|z_n|<R\exp\left({|z'|^2}\right)\Bigr\},
$$
for $0\le r<R\le\infty$. Every $CR$-isomorphism between two $R_{\delta}$-orbits is a composition of a map from of the form (\ref{gdelta}) and a dilation in the $z_n$-variable. Therefore, $F$ extends to a holomorphic automorphism of $\CC^n$. Hence, we obtain that $M$ is holomorphically equivalent to $D_r^R$ for some $0\le r<R\le\infty$ and therefore either to 
$$
D_{r/R,\,1}:=\Bigl\{(z',z_n)\in\CC^{n-1}\times\CC: r/R\exp\left({|z'|^2}\right)<|z_n|<
\exp\left({|z'|^2}\right)\Bigr\},
$$
or (for $R=\infty$) to 
$$
D_{0,-1}:=\Bigl\{(z',z_n)\in\CC^{n-1}\times\CC: 0<|z_n|<\exp\left({-|z'|^2}\right)\Bigr\}.
$$
This is, however, impossible since $d(D_{r/R,\,1})=d(D_{0,-1})=n^2$ (see e.g. \cite{I1}).

Assume now that ${\frak m}=\omega$. In this case ${\cal D}$ is the cylinder\linebreak ${\cal C}:=\left\{(z',z_n)\in\CC^{n-1}\times\CC:|z'|<1\right\}$. The $R_{\omega}$-orbit of every point in ${\cal C}$ has the form 
$$
\left\{(z',z_n)\in\CC^{n-1}\times\CC: |z'|^2+r\exp\left(\hbox{Re}\,z_n\right)=1\right\},
$$
where $r>0$, and any $R_{\omega}$-invariant domain in ${\cal C}$ is of the form
$$
\begin{array}{lll}
\Omega_r^R:&=&\Bigl\{(z',z_n)\in\CC^{n-1}\times\CC: |z'|<1,\\
&&r(1-|z'|^2)<\exp\left(\hbox{Re}\,z_n\right)<R(1-|z'|^2)\Bigr\},
\end{array}
$$
for $0\le r<R\le\infty$. Every $CR$-isomorphism between two $R_{\omega}$-orbits is a composition of a map from of the form (\ref{autgrpcov}) and a translation in the $z_n$-variable. Therefore, $F$ extends to a holomorphic automorphism of ${\cal C}$. In this case $M$ is holomorphically equivalent to $\Omega_r^R$ for some $0\le r<R\le\infty$, and hence either to
$$
\begin{array}{ll} 
\Omega_{r/R,1}:=&\Bigl\{(z',z_n)\in\CC^{n-1}\times\CC: |z'|<1,\,r/R(1-|z'|^2)<\\
&\exp\left(\hbox{Re}\,z_n\right)<(1-|z'|^2)\Bigr\},
\end{array}
$$
or  (for $R=\infty$) to
$$
\Omega_{0,-1}:=\Bigl\{(z',z_n)\in\CC^{n-1}\times\CC: |z'|<1,\,\exp\left(\hbox{Re}\,z_n\right)<(1-|z'|^2)^{-1}\Bigr\}.
$$
As before, this is impossible since $d(\Omega_{r/R,\,1})=d(\Omega_{0,-1})=n^2$ (see e.g. \cite{I1}).

Assume now that ${\frak m}=\varepsilon_{\alpha}$ for some $\alpha>0$. Here ${\cal D}$ is the domain ${\cal C}':={\cal C}\setminus\{z_n=0\}$. 
The $R_{\varepsilon_{\alpha}}$-orbit of every point in ${\cal C}'$ is of the form
$$
\left\{(z',z_n)\in\CC^{n-1}\times\CC: |z'|^2+r|z_n|^{\alpha}=1,\, z_n\ne 0\right\},
$$
where $r>0$, and every $R_{\varepsilon_{\alpha}}$-invariant domain in ${\cal C}'$ is given by
$$
\begin{array}{lll}
{\cal E}_{r,1/\alpha}^R:&=&\Bigl\{(z',z_n)\in\CC^{n-1}\times\CC: |z'|<1,\\
&&r(1-|z'|^2)^{1/\alpha}<|z_n|<R(1-|z'|^2)^{1/\alpha}\Bigr\},
\end{array}
$$
for $0\le r<R\le\infty$. Since every $CR$-isomorphism between $R_{\varepsilon_{\alpha}}$-orbits is a composition of a map of the form (\ref{autvarepsilon}) and a dilation in the $z_n$-variable, the map $F$ extends to an automorphism of ${\cal C}'$. Thus, we have shown that $M$ is holomorphically equivalent to ${\cal E}_{r,1/\alpha}^R$ for some $0\le r<R\le\infty$, and hence either to 
$$
\begin{array}{ll}
{\cal E}_{r/R,1/\alpha}:=&\Bigl\{(z',z_n)\in\CC^{n-1}\times\CC: |z'|<1,\, r(1-|z'|^2)^{1/\alpha}<\\
&\hspace{1.3cm}|z_n|<(1-|z'|^2)^{1/\alpha}\Bigr\},
\end{array}
$$
or (for $R=\infty$) to
$$
{\cal E}_{0,-1/\alpha}:=\Bigl\{(z',z_n)\in\CC^{n-1}\times\CC: |z'|<1,\, 0<
|z_n|<(1-|z'|^2)^{-1/\alpha}\Bigr\}.
$$
As above, this is impossible since $d({\cal E}_{r/R,1/\alpha})=d({\cal E}_{0,-1/\alpha})=n^2$ (see e.g. \cite{I1}).

Assume now that both real and complex hypersurface orbits are present in $M$. Since the action of $G(M)$ on $M$ is proper, it follows that the orbit space $M/G(M)$ is homeomorphic to one of the following: $\RR$, $S^1$, $[0,1]$, $[0,1)$ (see \cite{M}, \cite{B-B}, \cite{AA1}, \cite{AA2}), and thus there can be no more than two complex hypersurface orbits in $M$. It follows from (iii) of Proposition \ref{dim}  and Lemma \ref{groupsdeterm} that the model for every real hypersurface orbit is $\varepsilon_{\alpha}$ for some $\alpha>0$, $\alpha\in\QQ$. Let $M'$ be the manifold obtained from $M$ by removing all complex hypersurface orbits. It then follows from the above considerations that $M'$ is holomorphically equivalent to ${\cal E}_{r,1/\alpha}^R$ for some $0\le r<R\le\infty$.

Let $f:M'\ra {\cal E}_{r,1/\alpha}^R$ be a biholomorphic map satisfying 
(\ref{equivar}) for all $g\in G(M)$, $q\in M'$ and some isomorphism $\varphi: G(M)\ra R_{\varepsilon_{\alpha}}$. The group $R_{\varepsilon_{\alpha}}$ in fact acts on all of ${\cal C}$, and the orbit of any point in ${\cal C}$ with $z_n=0$ is the complex hypersurface        
$$
c_0:=\left\{(z',z_n)\in\CC^{n-1}\times\CC: |z'|<1,\,z_n=0\right\}.
$$
For a point $s\in{\cal C}$ denote by $J_s$ the isotropy subgroup of $s$ under the action of $R_{\varepsilon_{\alpha}}$. If $s_0\in c_0$ and $s_0=(z'_0,0)$, $J_{s_0}$ is isomorphic to $H^n_{k_1,k_2}$, where $k_1/k_2=2/\alpha n$ and consists of all maps of the form (\ref{repsilon}) for which the transformations in the $z'$-variables form the isotropy subgroup of the point $z_0'$ in $\hbox{Aut}(B^{n-1})$.     

Fix $s_0=(z_0',0)\in c_0$ and let
$$
N_{s_0}:=\left\{s\in {\cal E}_{r,1/\alpha}^R:J_s\subset J_{s_0}\right\}.
$$
We have
$$
N_{s_0}=\left\{(z',z_n)\in\CC^{n-1}\times\CC: z'=z'_0,\,r(1-|z_0'|^2)^{1/\alpha}<|z_n|<R(1-|z_0'|^2)^{1/\alpha}\right\}.
$$
Thus, $N_{s_0}$ is either an annulus (possibly, with an infinite outer radius) or a punctured disk. In particular, $N_{s_0}$ is a complex curve in ${\cal C}'$. Since $J_{s_0}$ is a maximal compact subgroup of $R_{\varepsilon_{\alpha}}$, $\varphi^{-1}(J_{s_0})$ is a maximal compact subgroup of $G(M)$. Let $O$ be a complex hypersurface orbit in $M$. For $q\in O$ the subgroup $I_q$ is compact and has dimension $(n-1)^2=\hbox{dim}\,J_{s_0}$. Therefore, $\varphi^{-1}(J_{s_0})$ is conjugate to $I_q$ for every $q\in O$ (in particular, $I_q$ is connected), and hence there exists $q_0\in O$ such that $\varphi^{-1}(J_{s_0})=I_{q_0}$. Since the isotropy subgroups in $R_{\varepsilon_{\alpha}}$ of distinct points in $c_0$ do not coincide, such a point $q_0$ is unique.

Let
$$
K_{q_0}:=\left\{q\in M': I_q\subset I_{q_0}\right\}.  
$$
Clearly, $K_{q_0}=f^{-1}(N_{s_0})$. Thus, $K_{q_0}$ is a $I_{q_0}$-invariant complex curve in $M'$ equivalent to either an annulus or a punctured disk. By Bochner's theorem there exist a local holomorphic 
change of coordinates $F$ near $q_0$ on $M$ that identifies an $I_{q_0}$-invariant neighborhood $U$ of $q_0$ with an $L_{q_0}$-invariant neighborhood of the origin in $T_{q_0}(M)$ such that $F(q_0)=0$ and $F(gq)=\alpha_{q_0}(g)F(q)$ for all $g\in I_{q_0}$ and $q\in U$ (here $L_{q_0}$ is the linear isotropy group and $\alpha_{q_0}$ is the isotropy representation at $q_0$). In the proof of Proposition \ref{dim} (see Case 1) we have seen that $L_{q_0}$ has two invariant subspaces in $T_{q_0}(M)$. One of them corresponds in our coordinates to $O$, the other to a complex curve $C$ intersecting $O$ at $q_0$. Observe that near $q_0$ the curve $C$ coincides with $K_{q_0}\cup\{q_0\}$. Therefore, in a neighborhood of $q_0$ the curve $K_{q_0}$ is a punctured analytic disk. Further, if a sequence $\{q_n\}$ from $K_{q_0}$ accumulates to $q_0$, the sequence $\{f(q_n)\}$ accumulates to one of the two ends of $N_{s_0}$, and therefore we have either $r=0$ or $R=\infty$. Since both these conditions cannot be satisfied simultaneously due to hyperbolicity of $M$, we conclude that $O$ is the only complex hypersurface orbit in $M$.

Assume first that $r=0$. We will extend $f$ to a map from $M$ onto the domain
\begin{equation}
\left\{(z',z_n)\in\CC^{n-1}\times\CC: |z'|^2+\frac{1}{R}|z_n|^{\alpha}<1\right\}\label{ellr}
\end{equation}
by setting $f(q_0)=s_0$, where $q_0\in O$ and $s_0\in c_0$ are related as specified above. The extended map is one-to-one and satisfies (\ref{equivar}) for all $g\in G(M)$, $q\in M$. To prove that $f$ is holomorphic on all of $M$, it suffices to show that $f$ is continuous on $O$. It will be more convenient for us to show that $f^{-1}$ is continuous on $c_0$. Let first $\{s_j\}$ be a sequence of points in $c_0$ converging to $s_0$. Then there exists a sequence $\{g_j\}$ of elements of $R_{\varepsilon_{\alpha}}$ converging to the identity such that $s_j=g_js_0$ for all $j$. Then $f^{-1}(s_j)=\varphi^{-1}(g_j)q_0$, and, since $\left\{\varphi^{-1}(g_j)\right\}$ converges to the identity, we obtain that $\{f^{-1}(s_j)\}$ converges to $q_0$. Next, let $\{s_j\}$ be a sequence of points in ${\cal E}_{0,\alpha}^R$ converging to $s_0$. Then we can find a sequence $\{g_j\}$ of elements of $R_{\varepsilon_{\alpha}}$ converging to the identity such that $g_js_j\in N_{s_0}$ for all $j$. Clearly, the sequence $\{f^{-1}(g_js_j)\}$ converges to $q_0$, and hence the sequence $\{f^{-1}(s_j)\}$ converges to $q_0$ as well. Thus, we have shown that $M$ is holomorphically equivalent to domain (\ref{ellr}) and hence to the domain
$$
 E_{\alpha}:=\left\{(z',z_n)\in\CC^{n-1}\times\CC: |z'|^2+|z_n|^{\alpha}<1\right\}.
$$
This is, however, impossible since $d(E_{\alpha})\ge n^2$.

Assume now that $R=\infty$. Observe that the action of the group $R_{\varepsilon_{\alpha}}$ on ${\cal C}$ extends to an action on $\tilde{\cal C}:=B^{n-1}\times\CC\PP^1$ by holomorphic transformations by setting $g(z',\infty):=(a(z'),\infty)$ for every $g\in R_{\varepsilon_{\alpha}}$, where $a$ is the corresponding automorphism of $B^{n-1}$ in the $z'$-variables (see formula (\ref{repsilon})). Now arguing as in the case $r=0$, we can extend $f$ to a biholomorphic map between $M$ and the domain in $\tilde{\cal C}$
$$
\left\{(z',z_n)\in\CC^{n-1}\times\CC: |z'|<1,\,|z_n|>r(1-|z'|^2)^{1/\alpha}\right\}\cup \left(B^{n-1}\times\{\infty\}\right).
$$
This domain is holomorphically equivalent to 
$$
{\cal E}_{-1/\alpha}:=\left\{(z',z_n)\in\CC^{n-1}\times\CC:|z'|<1,\, |z_n|<(1-|z'|^2)^{-1/\alpha}\right\},
$$
and so is $M$. This is, however, impossible since $d({\cal E}_{-1/\alpha})=n^2$. 

The proof of Proposition \ref{nospher} is complete.\qed

\section{The Case of Complex Hypersurface Orbits}\label{sectcomplex}
\setcounter{equation}{0}

We will now assume that all orbits in $M$ are complex hypersurfaces. As we have shown above, this is always the case for $n\ge 3$, unless $M$ is homogeneous. We will prove the following proposition.

\begin{proposition}\label{theoremcomplex}\sl Let $M$ be a connected hyperbolic manifold of dimension $n\ge 3$ with $d(M)=n^2-1$, and such that for every $p\in M$ its orbit $O(p)$ is a complex hypersurface in $M$. Then $M$ is holomorphically equivalent to $B^{n-1}\times S$, where $S$ is a hyperbolic Riemann surface with $d(S)=0$.
\end{proposition}

\noindent {\bf Proof:} Fix $p\in M$. It then follows from (iii) of Proposition \ref{dim} that $I_p^c$ is isomorphic to $U_{n-1}$, moreover, one can choose coordinates $(w_1,\dots,w_n)$ in $T_p(M)$ so that $L_p^c$ consists of all matrices of the form
\begin{equation}
\left(
\begin{array}{ll}
1 & 0\\
0 & B
\end{array}
\right),\label{formisotropy1}
\end{equation}
where $B\in U_{n-1}$ and $T_p(O(p))=\{w_1=0\}$. Arguing as in the proof of Lemma 4.4 of \cite{IKru1} we obtain that the full group $L_p$ consists of all matrices of the form
\begin{equation}
\left(
\begin{array}{ll}
\alpha & 0\\
0 & B
\end{array}
\right),\label{formisotropy}
\end{equation}
where $B\in U_{n-1}$ and $\alpha^m=1$ for some $m\ge 1$. It then follows (see e.g. Satz 4.3 of \cite{Ka}) that the kernel of the action of $G(M)$ on $O(p)$ is $J_p:=\alpha_p^{-1}(\ZZ_m)$, where we identify $\ZZ_m$ with the subgroup of $L_p$ that consists of all matrices of the form (\ref{formisotropy}) with $B=\hbox{id}$. Thus, $G(M)/J_p$ acts effectively on $O(p)$. Since $O(p)$ is holomorphically equivalent to $B^{n-1}$ and $\hbox{dim}\,G(M)=n^2-1=\hbox{dim}\,\hbox{Aut}(B^{n-1})$, we obtain that $G(M)/J_p$ is isomorphic to $\hbox{Aut}(B^{n-1})$. It then follows that $I_p$ is a maximal compact subgroup in $G(M)$ since its image under the projection $G(M)\ra\hbox{Aut}(B^{n-1})$ is a maximal compact subgroup of $\hbox{Aut}(B^{n-1})$. However, every maximal compact subgroup of a connected Lie group is connected whereas $I_p$ is not if $m>1$. Thus, $m=1$, hence $G(M)$ is isomorphic to $\hbox{Aut}(B^{n-1})$. In particular, $L_p$ fixes every point of the orthogonal complement $W_p$ to $T_p(O(p))$ in $T_p(M)$. Observe that the above arguments apply to every point in $M$.

Define
$$
N_p:=\left\{s\in M: I_s=I_p\right\}.
$$
Clearly, $I_p$ fixes every point in $N_p$ and $N_{gp}=gN_p$ for all $g\in G(M)$. Further, since for two distinct points $s_1,s_2$ lying in the same orbit we have $I_{s_1}\ne I_{s_2}$, the set $N_p$ intersects every orbit in $M$ at exactly one point. By Bochner's theorem there exist a local holomorphic 
change of coordinates $F$ near $p$ on $M$ that identifies an $I_p$-invariant neighborhood $U$ of $p$ with an $L_p$-invariant neighborhood $V$ of the origin in $T_p(M)$ such that $F(p)=0$ and $F(gq)=\alpha_p(g)F(q)$ for all $g\in I_p$ and $q\in U$. Since $L_p$ coincides with the group of matrices of the form (\ref{formisotropy1}), $N_p\cap U=F^{-1}(W_p\cap V)$. In particular, $N_p$ is a complex curve near $p$. Since the same argument can be carried out at every point of $N_p$, we obtain that $N_p$ is a closed complex hyperbolic curve in $M$.

We will now construct a biholomorphic map $\Phi: M\ra B^{n-1}\times N_p$. Let $\Psi: O(p)\ra B^{n-1}$ be a biholomorphism. For $q\in M$ let $r$ be the (unique) point where $N_p$ intersects $O(q)$. Let $g\in G(M)$ be such that $q=gr$. Then we set $\Phi(q):=(F(gp),r)$. By construction, $\Phi$ is biholomorphic. Since $M$ is holomorphically equivalent to $B^{n-1}\times N_p$, we have $d(N_p)=0$. 

The proof is complete.\qed

\section{The Homogeneous Case}\label{homog}
\setcounter{equation}{0}

In this section we will prove the following proposition.

\begin{proposition}\label{homogeneous} \sl If $M$ is a homogeneous connected hyperbolic manifold of dimension $n\ge 2$ with $d(M)=n^2-1$, then $n=4$ and $M$ is holomorphically equivalent to the tube domain
$$
\begin{array}{ll}
T=\Bigl\{(w_1,w_2,w_3,w_4)\in\CC^4:&(\hbox{Im}\,w_1)^2+(\hbox{Im}\,w_2)^2+
\\&(\hbox{Im}\,w_3)^2-(\hbox{Im}\,w_4)^2<0,\,\hbox{Im}\,w_4>0\Bigr\}.
\end{array}
$$
\end{proposition}

\noindent{\bf Proof:} The proof is similar to that of Proposition 5.1 of \cite{I1}. Since $M$ is homogeneous, by \cite{N}, \cite{P-S}, it is holomorphically equivalent to a Siegel domain $U$ of the second kind in $\CC^n$. For $n=2$, this gives that $M$ is equivalent to either $B^2$ or $\Delta^2$, which is impossible since $d(B^2)=8$ and $d(\Delta^2)=6$. For $n=3$ we obtain that $M$ is equivalent to one of the following domains: $B^3$, $B^2\times\Delta$, $\Delta^3$, $S$, where $S$ is the 3-dimensional Siegel space. None of these domains has an automorphism group of dimension 8.

Assume now that $n\ge 4$. The domain $U$ has the form
$$
U=\left\{(z,w)\in\CC^{n-k}\times\CC^k: \hbox{Im}\,w-F(z,z)\in C\right\},
$$
where $1\le k\le n$, $C$ is an open convex cone in $\RR^k$ not containing an entire affine line and $F=(F_1,\dots,F_k)$ is a $\CC^k$-valued Hermitian form on $\CC^{n-k}\times\CC^{n-k}$ such that $F(z,z)\in\overline{C}\setminus\{0\}$ for all non-zero $z\in\CC^{n-k}$.

We will show first that in most cases we have $k\le 2$. As we noted in \cite{IKra}
\begin{equation}
d(U)\le 4n-2k+\hbox{dim}\,{\frak g}_0(U).\label{form1}
\end{equation}
Here ${\frak g}_0(U)$ is the Lie algebra of all vector fields on $\CC^n$ of the form
$$
X_{A,B}=Az\frac{\partial}{\partial z}+Bw\frac{\partial}{\partial w},
$$
where $A\in{\frak{gl}}_{n-k}(\CC)$, $B$ belongs to the Lie algebra ${\frak g}(C)$ of the group $G(C)$ of linear automorphisms of the cone $C$, and the following holds
\begin{equation}
F(Az,z)+F(z,Az)=BF(z,z),\label{form3}
\end{equation}
for all $z\in\CC^{n-k}$. By the definition of Siegel domain, there exists a positive-definite linear combination $R$ of the components of the Hermitian form $F$. Then, for a fixed matrix $B$ in formula (\ref{form3}), the matrix $A$ is determined at most up to a matrix that is skew-Hermitian with respect to $R$. Since the dimension of the algebra of matrices skew-Hermitian with respect to $R$ is equal to $(n-k)^2$, we have
\begin{equation}
\hbox{dim}\,{\frak g}_0(U)\le (n-k)^2+\hbox{dim}\,{\frak g}(C).\label{estg}
\end{equation}

In Lemma 3.2 of \cite{IKra} we showed that 
\begin{equation}
\hbox{dim}\,{\frak g}(C)\le\frac{k^2}{2}-\frac{k}{2}+1.\label{lemmaestim}
\end{equation}
It now follows from (\ref{estg}) and (\ref{lemmaestim}) that the following holds      
$$
\hbox{dim}\,{\frak g}_0(U)\le\frac{3k^2}{2}-k\left(2n+\frac{1}{2}\right)+n^2+1,
$$
which together with (\ref{form1}) for gives
\begin{equation}
d(U)\le\frac{3k^2}{2}-k\left(2n+\frac{5}{2}\right)+n^2+4n+1.\label{form2}
\end{equation}
It is straightforward to check that the right-hand side of (\ref{form2}) is strictly less than $n^2-1$ if $k\ge 3$ for $n\ge 5$, and does not exceed 15 for $n=4$. Furthermore, for $n=4$ the right-hand side of (\ref{form2})  is equal to 15 only if $k=3$ or $k=4$ and $\hbox{dim}\,{\frak g}(C)=k^2/2-k/2+1$.

Suppose that $n=4$ and the right-hand side of (\ref{form2})  is equal to 15. In this case for every point $x_0\in C$ there exist coordinates in $\RR^k$ such that the isotropy subgroup of $x_0$ in $G(C)$ contains $SO_{k-1}(\RR)$ (see the proof of Lemma 3.2 in \cite{IKra}). Then after a linear change of coordinates the cone $C$ takes the form
$$
\left\{x=(x_1,\dots,x_k)\in\RR^k:\left\langle x,x\right\rangle <0,\,x_k>0\right\},
$$
where $\left\langle x,x\right\rangle:=x_1^2+\dots+x_{k-1}^2-x_k^2$. 
In these coordinates the algebra ${\frak g}(C)$ is generated by       
the subalgebra of scalar matrices in ${\frak {gl}}_k(\RR)$ and the algebra of pseudo-orthogonal matrices ${\frak o}_{k-1,1}(\RR)$. Assume first that $k=3$. Then we have $F=(v_1|z|^2,v_2|z|^2,v_3|z|^2)$ for some vector  $v:=(v_1,v_2,v_3)\in C$. It follows from (\ref{form3}) that $v$ is an eigenvector of the matrix $B$ for every $X_{A,B}\in{\frak g}_0(U)$, which implies that $\hbox{dim}\,{\frak g}_0(U)=3$. Hence by (\ref{form1}) we have $\hbox{dim}\,\hbox{Aut}(U)\le 13$, which is impossible.

Suppose now that $k=4$. In this case $U$ is holomorphically equivalent to the tube domain $T$. Let ${\frak g}(T)$ be the Lie algebra of $\hbox{Aut}(T)$. It follows from the results of \cite{KMO} that ${\frak g}(T)$ is a graded Lie algebra
$$
{\frak g}(T)={\frak g}_{-1}(T)\oplus{\frak g}_0(T)\oplus{\frak g}_1(T),
$$
where ${\frak g}_{-1}$ is spanned by $i\partial/\partial w_j$, $j=1,2,3,4$, and $\hbox{dim}\,{\frak g}_1(T)\le 4$. Clearly, ${\frak g}_0(T)$ is isomorphic to $\RR\oplus{\frak o}_{3,1}(\RR)$ and thus has dimension 7. The component ${\frak g}_1(T)$ also admits an explicit description (see e.g. p. 218 in \cite{S}). It follows from this description that ${\frak g}_1(T)$ consists of all vector fields of the form
$$
\begin{array}{ll}
Z_{\alpha,\beta,\gamma,\delta}:=&\displaystyle\Bigl(\alpha(w_1^2-w_2^2-w_3^2+w_4^2)+2(\beta w_1w_2+\gamma w_1w_3 +\delta w_1w_4)\Bigr)\frac{\partial}{\partial w_1}+\\
\vspace{0cm}&\\
&\displaystyle\Bigl(\beta(-w_1^2+w_2^2-w_3^2+w_4^2)+2(\alpha w_1w_2+\gamma w_2w_3 +\delta w_2w_4)\Bigr)\frac{\partial}{\partial w_2}+\\
\vspace{0cm}&\\
&\displaystyle\Bigl(\gamma(-w_1^2-w_2^2+w_3^2+w_4^2)+2(\alpha w_1w_3+\beta w_2w_3 +\delta w_3w_4)\Bigr)\frac{\partial}{\partial w_3}+\\
\vspace{0cm}&\\
&\displaystyle\Bigl(\delta(w_1^2+w_2^2+w_3^2+w_4^2)+2(\alpha w_1w_4+\beta w_2w_4 +\gamma w_3w_4)\Bigr)\frac{\partial}{\partial w_4},
\end{array}
$$
where $\alpha,\beta,\gamma,\delta\in\RR$, and thus has dimension 4. Therefore, $\hbox{dim}\,\hbox{Aut}(T)=15$. It is also clear that $T$ is homogeneous under affine automorphisms.      

Assume now that $n\ge 4$ is arbitrary and $k\le 2$. If $k=1$, the domain $U$ is equivalent to $B^n$ which is impossible. Hence $k=2$. It follows from (\ref{form3}) that the matrix $A$ is determined by the matrix $B$ up to a matrix $L\in{\frak {gl}}_{n-2}(\CC)$ satisfying
$$
F(Lz,z)+F(z,Lz)=0,
$$
for all $z\in\CC^{n-2}$. Let $s$ be the dimension of the subspace of all such matrices $L$. Then
$$
\hbox{dim}\,{\frak g}_0(U)\le s+\hbox{dim}\,{\frak g}(C),
$$
and (\ref{lemmaestim}) yields
$$
\hbox{dim}\,{\frak g}_0(U)\le s+2,
$$
which together with (\ref{form1}) implies
\begin{equation}
s\ge n^2-4n+1.\label{form4}
\end{equation}
 
By the definition of Siegel domain, there exists a positive-definite linear combination of the components of $F$, and we can assume that $F_1$ is positive-definite. Further, applying an appropriate linear transformation of the $z$-variables, we can assume that $F_1$ is given by the identity matrix and $F_2$ by a diagonal matrix.

Suppose first that the matrix of $F_2$ is scalar. If $F_2\equiv 0$, then $U$ is holomorphically equivalent to $B^{n-1}\times\Delta$ which is impossible. If $F_2\not\equiv 0$, then $U$ is holomorphically equivalent to the domain
$$
V:=\left\{(z,w)\in\CC^{n-2}\times\CC^2:\hbox{Im}\, w_1-|z|^2>0,\,\hbox{Im}\,w_2-|z|^2>0\right\}.
$$
It was shown in \cite{IKra} that $d(V)\le n^2-2n+3$ and hence $d(V)<n^2-1$. Thus, the matrix of $F_2$ is not scalar. Inequality (\ref{form4}) now yields that the matrix of $F_2$ can have at most one pair of distinct eigenvalues, and therefore $n=4$ and $U$ is holomorphically equivalent to $B^2\times B^2$. This is clearly impossible, and the proof of the proposition is complete.\qed

\section{Examples for the Case $n=2$, $d(M)=3$}\label{examples23}
\setcounter{equation}{0}

In this section we give examples of families of hyperbolic domains in $\CC^2$ and $\CC\PP^2$ with automorphism groups of dimension 3 whose orbit structure is different from that observed above for $n\ge 3$. Define
$$
\Omega_t:=\left\{(z,w)\in\CC^2: |z|^2+|w|^2-1< t |z^2+w^2-1|\right\},
$$
where $0<t\le 1$. Clearly, $\Omega_t$ is bounded if $0<t<1$. Further, $\Omega_1$ is hyperbolic since it is contained in the hyperbolic product domain
$$
\left\{(z,w)\in\CC^2: z,w\not\in(-\infty,-1]\cup[1,\infty)\right\}.
$$
The group $\hbox{Aut}(\Omega_t)$ for every $t$ consists of the maps
$$
\left(
\begin{array}{c}
z\\
w
\end{array}
\right) \mapsto \displaystyle\frac{\left(\begin{array}{cc}
a_{11} & a_{12}\\
a_{21} & a_{22}
\end{array}
\right)\left(
\begin{array}{c}
z\\
w
\end{array}\right)+\left(
\begin{array}{c}
b_1\\
b_2
\end{array}
\right)}{c_1z+c_2w+d},
$$
where
\begin{equation}
Q:=\left(\begin{array}{ccc}
a_{11} & a_{12} & b_1\\
a_{21} & a_{22} & b_2\\
c_1& c_2 &d
\end{array}\label{matq}
\right)
\in SO_{2,1}(\RR),
\end{equation}
and thus is 3-dimensional. The group $\hbox{Aut}(\Omega_t)$ has two connected components (that correspond to the connected components of $SO_{2,1}(\RR)$), and its identity component $G(\Omega_t)$ is given by the condition $a_{11}a_{22}-a_{12}a_{21}>0$. The orbits of $G(\Omega_{t})$ on $\Omega_t$ are as follows:
$$
\begin{array}{ll}
O^{\Omega}_{\alpha}:=&\left\{(z,w)\in\CC^2: |z|^2+|w|^2-1=\alpha|z^2+w^2-1|\right\}\setminus\\
&\left\{(x,u)\in\RR^2:x^2+u^2=1\right\},\quad -1<\alpha<t,\\
\vspace{0cm}&\\
\Delta_{\RR}:=&\left\{(x,u)\in\RR^2:x^2+u^2<1\right\}.
\end{array}
$$
Note that $O^{\Omega}_0$ is the only spherical real hypersurface orbit in $\Omega_t$ and that $\Delta_{\RR}$ is a totally real orbit. All the orbits are pairwise $CR$ non-equivalent.

The next family of domains is associated with a different action of $SO_{2,1}(\RR)$ on a part of $\CC^2$. Define
$$
D_t:=\left\{(z,w)\in\CC^2: 1+|z|^2-|w|^2> t |1+z^2-w^2|,\,\hbox{Im}\,z(1+\overline{w})>0\right\},
$$
where $t\ge 1$. All these domains lie in the hyperbolic product domain 
$$
\left\{(z,w)\in\CC^2: \hbox{Im}\,z>0,\,w\not\in(-\infty,-1]\cup[1,\infty)\right\},
$$
hence they are hyperbolic as well. For every matrix $Q\in SO^c_{2,1}(\RR)$ as in (\ref{matq}) consider the map
$$
\left(
\begin{array}{c}
z\\
w
\end{array}
\right) \mapsto \displaystyle\frac{\left(\begin{array}{cc}
a_{22} & b_2\\
c_2 & d
\end{array}
\right)\left(
\begin{array}{c}
z\\
w
\end{array}\right)+\left(
\begin{array}{c}
a_{21}\\
c_1
\end{array}
\right)}{a_{12}z+b_1w+a_{11}}.
$$
The group $\hbox{Aut}(D_t)=G(D_t)$ for every $t$ consists of all such maps. The orbits of $G(D_t)$ on $D_t$ are the following non-spherical hypersurfaces
$$
\begin{array}{ll}
O^D_{\alpha}:=&\Bigl\{(z,w)\in\CC^2: 1+|z|^2-|w|^2=\alpha|1+z^2-w^2|,\\
&\hbox{Im}\,z(1+\overline{w})>0\Bigr\},\quad \alpha>t.
\end{array}
$$
All the orbits are pairwise $CR$ non-equivalent.

The next family of domains is associated with an action of $SO_3(\RR)$ on $\CC\PP^2$. Define
$$
E_t:=\left\{(z:w:\zeta)\in\CC\PP^2: |z|^2+|w|^2+|\zeta|^2< t |z^2+w^2+\zeta^2|\right\},
$$ 
where $t>1$. The domain $E_t$ is hyperbolic for each $t$ since it is covered in a 2-to-1 fashion by the manifold
$$
\Bigl\{(z,w,\zeta)\in\CC^3: |z|^2+|w|^2+|\zeta|^2<t,\,z^2+w^2+\zeta^2=1\Bigr\},
$$
which is clearly hyperbolic; the covering map is $(z,w,\zeta)\mapsto (z:w:\zeta)$. The group $\hbox{Aut}(E_t)=G(E_t)$ for every $t$ is given by applying matrices from $SO_3(\RR)$ to vectors of homogeneous coordinates. The action of the group $G(E_t)$ on $E_t$ has the totally real orbit $\RR\PP^2$, and the rest of the orbits are the following non-spherical hypersurfaces
$$
O^E_{\alpha}:=\left\{(z:w:\zeta)\in\CC\PP^2: |z|^2+|w|^2+|\zeta|^2=\alpha |z^2+w^2+\zeta^2|\right\}, \quad 1<\alpha<t.
$$
All the orbits are pairwise $CR$ non-equivalent.

Next, define
$$
S_t:=\left\{(z,w)\in\CC^2: \left(\hbox{Re}\,z\right)^2+\left(\hbox{Re}\,w\right)^2<t\right\},
$$
where $t>0$. All these domains are clearly hyperbolic and the group $\hbox{Aut}(S_t)$ for every $t$ consists of all maps of the form
$$
\left(
\begin{array}{l}
z\\
w
\end{array}
\right)\mapsto C
\left(
\begin{array}{l}
z\\
w
\end{array}
\right)+i
\left(
\begin{array}{l}
p\\
q
\end{array}
\right),
$$
where $C\in O_2(\RR)$ and $p,q\in\RR$. The group $G(S_t)$ is given by matrices $C\in SO_2(\RR)$. The action of the group $G(S_t)$ on $S_t$ has the totally real orbit
$$
\left\{(z,w)\in\CC^2: \hbox{Re}\,z=0,\,\hbox{Re}\,w=0\right\}, 
$$
and the rest of the orbits are the following non-spherical tube hypersurfaces
$$
O^S_{\alpha}:=\left\{(z,w)\in\CC^2: \left(\hbox{Re}\,z\right)^2+\left(\hbox{Re}\,w\right)^2=\alpha\right\},\quad 0<\alpha<t.
$$
Every non-spherical orbit is clearly $CR$-equivalent to $O^S_1$.

Now fix $a\in\RR$ such that $|a|>1$, $a\ne 1,2$, and consider the following family of tube domains
$$
R_{a,t}:=\left\{(z,w)\in\CC^2: \hbox{Re}\,z<t\left(\hbox{Re}\,w\right)^a,\,\hbox{Re}\,w>0\right\},
$$
where $t>0$. All these domains are obviously hyperbolic and the group $\hbox{Aut}(R_{a,t})=G(R_{a,t})$ consists of all the maps
$$
\left(
\begin{array}{l}
z\\
w
\end{array}
\right)\mapsto
\left(
\begin{array}{l}
\lambda^a z\\
\lambda w
\end{array}
\right)+i
\left(
\begin{array}{l}
p\\
q
\end{array}
\right),
$$
where $\lambda>0$ and $p,q\in\RR$. The action of this group on $R_{a,t}$ has the Levi-flat orbit
$$
\left\{(z,w)\in\CC^2: \hbox{Re}\,z=0,\,\hbox{Re}\,w>0\right\}, 
$$
which is foliated by the half-planes
$$
\left\{(z,w)\in\CC^2: z=ic,\,\hbox{Re}\,w>0\right\}, \quad c\in\RR. 
$$
All other orbits are the following non-spherical hypersurfaces
$$
O^R_{a,\alpha}:=\left\{(z,w)\in\CC^2: \hbox{Re}\,z=\alpha\left(\hbox{Re}\,w\right)^a,\,\hbox{Re}\,w>0\right\},\quad\alpha<t,\,\alpha\ne 0.
$$
Every non-spherical orbit is $CR$-equivalent to $O^R_{a,1}$.

Further, define
$$
U_t:=\left\{(z,w)\in\CC^2:\hbox{Re}\,z<\hbox{Re}\,w\cdot\ln\left(t\hbox{Re}\,w\right),\,\hbox{Re}\,w>0\right\},
$$
where $t>0$. All these domains are clearly hyperbolic and the group $\hbox{Aut}(U_t)=G(U_t)$ consists of all the maps
$$
\left(
\begin{array}{l}
z\\
w
\end{array}
\right)\mapsto
\left(
\begin{array}{l}
\lambda z+(\lambda\ln\lambda)w\\
\lambda w
\end{array}
\right)+i
\left(
\begin{array}{l}
p\\
q
\end{array}
\right),
$$
where $\lambda>0$ and $p,q\in\RR$. The orbits of $G(U_t)$ on $U_t$ are the following non-spherical hypersurfaces
$$
O^U_{\alpha}:=\left\{(z,w)\in\CC^2:\hbox{Re}\,z=\hbox{Re}\,w\cdot\ln\left(\alpha\hbox{Re}\,w\right),\,\hbox{Re}\,w>0\right\},\quad 0<\alpha<t.
$$
Every orbit is $CR$-equivalent to $O^U_1$.

Finally, fix $a>0$ and consider
$$
V_{a,t,s}:=\left\{(z,w)\in\CC^2:se^{a\varphi}<r<te^{a\varphi}\right\},
$$   
where $t>0$, $e^{-2\pi a}t<s<t$, and $(r,\varphi)$ denote the polar coordinates in the $(\hbox{Re}\,z,\hbox{Re}\,w)$-plane with $\varphi$ varying from $-\infty$ to $\infty$ (thus, the boundary of $V_{a,t,s}$ consists of two infinite spirals). All these domains are hyperbolic and $\hbox{Aut}(V_{a,t,s})=G(V_{a,t,s})$ consists of all maps of the form
$$
\left(
\begin{array}{l}
z\\
w
\end{array}
\right)\mapsto e^{a\beta}
\left(
\begin{array}{rr}
\cos\beta & \sin\beta\\
-\sin\beta & \cos\beta
\end{array}
\right)
\left(
\begin{array}{l}
z\\
w
\end{array}
\right)+i
\left(
\begin{array}{l}
p\\
q
\end{array}
\right),
$$
where $\beta,p,q\in\RR$. The orbits under the action of $G(V_{a,t,s})$ on $V_{a,t,s}$ are the following non-spherical hypersurfaces
$$
O^V_{a,\alpha}:=\left\{(z,w)\in\CC^2:r=\alpha e^{a\varphi}\right\},\quad s<\alpha<t.
$$
Clearly, every orbit is $CR$-equivalent to $O^V_{a,1}$. 

The orbits $O^{\Omega}_{\alpha}$ with $-1<\alpha<1$ and $\alpha\ne 0$, $O^D_{\alpha}$ with $\alpha>1$, $O^E_{\alpha}$ with $\alpha>1$, $O^S_1$, $O^R_{a,1}$ with $|a|>1$ and $a\ne 1,2$, $O^U_1$, $O^V_{a,1}$ with $a>0$ are part of E. Cartan's classification of homogeneous hypersurfaces in the non-spherical case (see \cite{C}). They are pairwise $CR$ non-equivalent, both locally and globally, and give a complete classification from the local point of view. To obtain a global classification, one has to additionally consider all possible covers of these hypersurfaces (see \cite{I2}).

We will now give an example of a hyperbolic domain in $\CC^2$, for which almost every orbit is spherical. Define
$$
W_t:=\Bigl\{(z,w)\in\CC^2: \hbox{Re}\,w>|z|^2+t\left(\hbox{Re}\,z\right)^2,\,\hbox{Re}\,z>0\Bigr\},
$$
where $t\in\RR$. This domain is hyperbolic since for $t\ge -2$ it is equivalent to a subdomain of the hyperbolic product domain
$$
\left\{(z,w)\in\CC^2: \hbox{Re}\,z>0,\,\hbox{Re}\,w>0\right\},
$$
and for $t<-2$ it is equivalent to the hyperbolic domain
$$
\left\{(z,w)\in\CC^2: \hbox{Re}\,w<|z|^2,\,\hbox{Re}\,z>0\right\}.
$$ 
The group $\hbox{Aut}(W_t)=G(W_t)$ consists of the maps
$$
\begin{array}{lll}
z & \mapsto & \lambda z+ia,\\
w & \mapsto & \lambda^2w-2i\lambda a z+a^2+i\beta,
\end{array}
$$
where $\lambda>0$, $a,\beta\in\RR$ (cf. (\ref{thegroupsphpt})). The action of this group on $W_t$ has the Levi-flat orbit
$$
\left\{(z,w)\in\CC^2: \hbox{Re}\,(w+z^2)=0,\,\hbox{Re}\,z>0\right\}, 
$$
which is foliated by the complex curves
$$
\left\{(z,w)\in\CC^2: w+z^2=ic,\,\hbox{Re}\,z>0\right\}, \quad c\in\RR. 
$$
All other orbits are the following spherical hypersurfaces
$$
O^W_{\alpha}:=\Bigl\{(z,w)\in\CC^2:\hbox{Re}\,(w-\alpha z^2/2 )=|z|^2(1+\alpha/2),\,\hbox{Re}\,z>0\Bigr\},\quad \alpha>t. 
$$
Clearly, every spherical orbit is $CR$-equivalent to $O^W_0$.

{\obeylines
Department of Mathematics
The Australian National University
Canberra, ACT 0200
AUSTRALIA
E-mail: alexander.isaev@maths.anu.edu.au
}

\end{document}